\documentclass[12pt]{amsart}
\usepackage[active]{srcltx}
\usepackage{amssymb,amsmath}
\usepackage{graphicx}
\usepackage{amsbsy}
\usepackage{mathrsfs}

\newtheorem{theorem}{Theorem}[section]

\newtheorem{remark}{Remark}[section]
\numberwithin{equation}{section}

\begin{document}
\begin{center}
  {\Large\bf
 Elliptic problems on the ball endowed with Funk-type metrics}\\
\vspace{0.5cm} {  Alexandru Krist\'aly$^{*,**}$\ \&\ Imre J. Rudas$^{**}$\\
\vspace{0.3cm} {\small  $^{*}$Department of Economics, Babe\c
s-Bolyai University, Cluj-Napoca, Romania\\Email address:
alexandrukristaly@yahoo.com\\ $^{**}$Institute of Applied
Mathematics, \'Obuda University, 1034 Budapest, Hungary}}
\end{center}

\vspace{1cm}

\noindent {\bf Abstract.} {\footnotesize \noindent We study Sobolev
spaces on the $n-$dimensional unit ball $B^n(1)$ endowed with a
parameter-depending Finsler metric $F_a$, $a\in [0,1],$ which
interpolates between the Klein metric $(a=0)$ and Funk metric
$(a=1)$, respectively. We show that the standard Sobolev space
defined on the Finsler manifold $(B^n(1),F_a)$ is a vector space if
and only if $a\in [0,1).$ Furthermore, by exploiting variational
arguments, we provide non-existence and existence results for
sublinear elliptic problems on $(B^n(1),F_a)$ involving the
Finsler-Laplace operator whenever $a\in [0,1).$ }

\vspace{1cm}

%\noindent {\bf MSC}: 58J05, 53C60, 58J60.

\section{Introduction}

The theory of Sobolev spaces on complete Riemannian manifolds is
well understood and widely applied into the study of various
elliptic problems, see e.g. Hebey \cite{Hebey}. Although Finsler
geometry is a natural extension of Riemannian geometry, very little
is known about Sobolev spaces on non-compact Finsler manifolds. One
of the most important features of Finsler structures is that they
can describe non-reversible phenomena. Such examples appear for
instance in the Matsumoto mountain slope metric (describing the law
of walking on a mountain under the action of gravity), the
Poincar\'e-Finsler disc model, the Funk metric on the
$n-$dimensional unit ball $B^n(1)\subset \mathbb R^n$, etc.

The purpose of this paper is to study properties of the Sobolev
space and some elliptic problems involving the Finsler-Laplace
operator on $B^n(1)\subset \mathbb R^n$ which is endowed with a
Funk-type metric. To be more precise, let $B^n(1)=\{x\in \mathbb
R^n: |x|<1\}$ be the $n-$dimensional unit ball, $n\geq 2$, and for
every $a\in [0,1]$, we consider the function $F_a:B^n(1)\times
\mathbb R^{n}\to \mathbb R$ defined by
\begin{equation}\label{F-a-metrika}
    F_a(x,y)=\frac{\sqrt{|y|^2-(|x|^2|y|^2-\langle
x,y\rangle^2)}}{1-|x|^2}+a\frac{\langle x,y\rangle}{1-|x|^2},\ x\in
B^n(1),\ y\in T_xB^n(1)=\mathbb R^n.
\end{equation}
 Hereafter, $|\cdot|$ and
$\langle\cdot, \cdot\rangle$ denote the $n-$dimensional Euclidean
norm and inner product.  Standard arguments from Finsler geometry
show that the pair $(B^n(1),F_a)$ is a Finsler manifold of Randers
type. In fact, for $a=0$, the manifold $(B^n(1),F_0)$ reduces to the
well known Riemannian {\it Klein model}; for $a=1$, $(B^n(1),F_1)$
is the usual Finslerian {\it Funk model}, see Cheng and Shen
\cite{Cheng-Shen} and Shen \cite{Shen}. We introduce the Sobolev
space associated with
  $(B^n(1),F_a)$; namely, let
%$$W^{1,2,a}(B^n(1))=\left\{u\in W^{1,2}_{\rm loc}(B^n(1)):\displaystyle
%\int_{B^n(1)} F_a^{*2}(x,Du(x)){\rm d}V_{F_a}(x)<+\infty\right\},$$
 $W_{0}^{1,2,a}(B^n(1))$ be the closure of $C_0^\infty(B^n(1))$
%in $W^{1,2,a}(B^n(1))$
with respect to the (positively homogeneous)
norm
\begin{equation}\label{Sobolev-norm}
    \|u\|_{W^{1,2,a}}=\left(\displaystyle \int_{B^n(1)} F_a^{*2}(x,Du(x)){\rm
d}V_{F_a}(x)+\displaystyle \int_{B^n(1)} u^2(x){\rm
d}V_{F_a}(x)\right)^{1/2},
\end{equation}
where $F_a^*$ denotes the polar transform  of $F_a,$ and ${\rm
d}V_{F_a}$ is the canonical Hausdorff measure on $(B^n(1),F_a)$; see
Section \ref{sect-2} for details. For $a=0$ (Klein metric case), the
space $W_{0}^{1,2,0}(B^n(1))$ endowed with the (absolutely
homogenous) norm $\|\cdot\|_{H_1^2}=\|\cdot\|_{W^{1,2,0}}$ is
nothing but the usual Sobolev space $H_1^2(B^n(1))$
 on the Riemannian manifold $(B^n(1),F_0)$ having a Hilbert structure, see Hebey
\cite{Hebey}. Sobolev spaces on generic Finsler manifolds were
introduced (and studied in the compact case) in \cite{Ge-Shen} and
  \cite{Ohta-Sturm}.

Our first result reads as follows.

\begin{theorem}\label{theorem-nem-vektor-ter}
 Let $a\in [0,1]$. Then the following assertions are equivalent:
\begin{itemize}
  \item[(i)] $W_{0}^{1,2,a}(B^n(1))$ is a vector space over $\mathbb R;$
  \item[(ii)] $a\in [0,1).$
\end{itemize}
\end{theorem}
\noindent The proof of Theorem \ref{theorem-nem-vektor-ter} is based
on the following two facts:
\begin{itemize}
  \item when $a=1,$ i.e., $F_a$ is the usual Funk  metric, we
  construct a function  $u\in W_{0}^{1,2,1}(B^n(1))$ such that $-u\notin
  W_{0}^{1,2,1}(B^n(1))$;
  \item when $a\in [0,1)$, the vector space structure of
  $W_{0}^{1,2,a}(B^n(1))$ follows in a standard way by exploiting the convexity of $F_a^*$ and the finiteness of the reversibility constant of $(B^n(1),F_a)$. As a byproduct,
   the equivalence of the norms $\|\cdot\|_{W^{1,2,a}}$ and
   $\|\cdot\|_{H_1^2}$ easily follows.
\end{itemize}

In the second part of the paper we consider the highly nonlinear
problem
\[ \   \left\{ \begin{array}{lll}
 -\boldsymbol{\Delta}_{F_a} u=\lambda \kappa(x) g(u) &\mbox{in} &  B^n(1); \\
% u\geq 0 &\mbox{in} &   \Omega;\\
u\to 0 &\mbox{if} &  |x|\to 1,
 \end{array}\right. \eqno{({\mathcal P}_{\lambda})}\]
where $a\in [0,1)$, $\boldsymbol{\Delta}_{F_a}$ denotes the
Finsler-Laplace operator
 associated with the Funk-type metric $F_a$ on $B^n(1)$, $\lambda\geq 0$ is a
 parameter, $\kappa\in L^1(B^n(1))\cap L^\infty(B^n(1))$
 and $g:[0,\infty)\to \mathbb R$ is a continuous function.
 Note that when $a=0$,
$\boldsymbol{\Delta}_{F_a}$ becomes the usual Laplace-Betrami
operator $\Delta_{F_0}$ on the Klein ball model  $(B^n(1),F_0)$;
however, when $a\neq 0$, the operator $\boldsymbol{\Delta}_{F_a}$ is
highly non-linear (neither additive nor absolutely homogeneous). On
the continuous function $g$, we require
\begin{itemize}
   \item[(g1)] $g(s)=o(s)$ as $s\to 0^+$ and $s\to \infty;$
   \item[(g2)] $G(s_0)>0$ for some $s_0>0,$ where $G(s)=\int_0^s g(t)dt.$
 \end{itemize}
 The assumption $g(s)=o(s)$ as $s\to
\infty$ means that $g$ is sublinear at infinity.  Moreover, due to
{\rm (g1)} and {\rm
  (g2)},
the number $$c_g=\max_{s> 0}\frac{g(s)}{s}$$ is well-defined and
positive.

Our second result reads as follows.

\begin{theorem}\label{theorem-sublinear}   Let $a\in [0,1)$,
$\kappa\in L^1(B^n(1))\cap L^\infty(B^n(1))\setminus\{0\}$ be a
radially symmetric non-negative function,
 and a continuous
  function $g:[0,\infty)\to \mathbb R$ verifying {\rm (g1)} and {\rm
  (g2)}.
 Then
\begin{itemize}
  \item[{\rm (i)}]   $(\mathcal P_{\lambda})$ has only the zero
solution whenever  $0\leq
\lambda<c_g^{-1}\|\kappa\|_{L^\infty}^{-1}\frac{(n-1)^2(1-a^2)^\frac{n+1}{2}}{4(1+a)^2};$

  \item[{\rm (ii)}]  there exists $\tilde \lambda>0$ such that  $(\mathcal P_{\lambda})$ has at least
two distinct  non-zero, non-negative, radially symmetric weak
solutions whenever
  $\lambda>\tilde\lambda$.
\end{itemize}
\end{theorem}

The proof of (i) combines a direct computation with a result of
Federer and Fleming \cite{FF} applied for the Klein ball model. In
order to prove (ii), we shall exploit Theorem
\ref{theorem-nem-vektor-ter} together with variational arguments
(minimization, mountain pass and the principle of symmetric
criticality). In fact, Theorem \ref{theorem-sublinear} seems to be
the first existence result within the class of elliptic problems on
{\it non}-compact Finsler manifolds.

\section{Preliminaries}\label{sect-2}
\subsection{Randers spaces.}
Let $M$ be a smooth, $n-$dimensional manifold and $TM=\bigcup_{x \in
M}T_{x} M $ be its tangent bundle. Throughout of this subsection,
the function $F:TM\to [0,\infty)$ is given by
\begin{equation}\label{Randers-metrika}
    F(x,y)=\sqrt{h_x(y,y)}+\beta_x(y),\ (x,y)\in TM,
\end{equation}
where $h$ is a Riemannian metric on $M$, $\beta$ is an 1-form on
$M$, and we assume that
$$\|\beta\|_h(x)=\sqrt{h_x^*(\beta_x,\beta_x)}<1,\ \forall x\in M.$$
Here, the co-metric $h^*_x$ can be identified by  the inverse of the
symmetric, positive definite matrix $h_x$. The pair $(M,F)$ is a
{\it Randers space} which is a typical Finsler manifold, i.e., the
following properties hold:

(a) $F\in C^{\infty}(TM\setminus\{ 0 \});$

(b) $F(x,ty)=tF(x,y)$ for all $t\geq 0$ and $(x,y)\in TM;$

(c) $g_{(x,y)}=[g_{ij}(x,y)]:=\left[\frac12F^{2}%
(x,y)\right]_{y^{i}y^{j}}$ is positive definite for all $(x,y)\in
TM\setminus\{ 0 \},$\\
\noindent see Bao, Chern and Shen \cite{BCS}. Clearly, the Randers
metric $F$ in (\ref{Randers-metrika}) is symmetric, i.e.,
$F(x,-y)=F(x,y)$ for every $(x,y)\in TM,$ if and only if $\beta=0$.

Unlike the Levi-Civita connection \index{Levi-Civita connection} in
the Riemannian case, there is no unique natural connection in the
Finsler geometry. Among these connections on the pull-back bundle
$\pi ^{*}TM,$ we choose a torsion free and almost metric-compatible
linear connection on $\pi ^{*}TM$, the so-called \textit{Chern
connection}, see Bao, Chern and Shen \cite[Theorem 2.4.1]{BCS}.
Since the notions of geodesics, forward/backward completeness and
flag curvature will not be explicitly used in the sequel, we assume
the reader is familiar with them; for details, see \cite{BCS}.

Let $\sigma: [0,r]\to M$ be a piecewise $C^{\infty}$ curve. The
value $ L_F(\sigma)= \displaystyle\int_{0}^{r} F(\sigma(t),
\dot\sigma(t))\,{\text d}t $ denotes the \textit{integral length} of
$\sigma.$  For $x_1,x_2\in M$, denote by $\Lambda(x_1,x_2)$ the set
of all piecewise $C^{\infty}$ curves $\sigma:[0,r]\to M$ such that
$\sigma(0)=x_1$ and $\sigma(r)=x_2$. Define the {\it distance
function} $d_{F}: M\times M \to[0,\infty)$ by
\begin{equation}\label{quasi-metric}
  d_{F}(x_1,x_2) = \inf_{\sigma\in\Lambda(x_1,x_2)}
 L_F(\sigma).
\end{equation}
One clearly has that $d_{F}(x_1,x_2) =0$ if and only if $x_1=x_2,$
 and $d_F$ verifies the triangle inequality.

  The {\it Hausdorff  volume form} ${\text
d}V_F$ on the Randers space $(M,F)$ is given by
\begin{equation}\label{randers-volume}
{\text d}V_F(x)=\left(1-\|\beta\|^2_h(x)\right)^\frac{n+1}{2}{\text
d}V_h(x),
\end{equation}
where ${\text d}V_h(x)$ denotes the usual Riemannian volume form of
$h$ on $M,$ see Cheng and Shen \cite{Cheng-Shen}.

For every $(x,\alpha)\in T^*M$, the {\it polar transform} (or,
co-metric) of $F$ from (\ref{Randers-metrika}) is
\begin{eqnarray}\label{polar-transform}
 \nonumber F^*(x,\alpha)&\stackrel{\rm def.}{=}& \sup_{y\in T_xM\setminus
    \{0\}}\frac{\alpha(y)}{F(x,y)} \\
 &=&\frac{\sqrt{h_x^{*2}(\alpha,\beta)+(1-\|\beta\|_h^2(x))\|\alpha\|_h^2(x)}-h_x^{*}(\alpha,\beta)}{1-\|\beta\|_h^2(x)}.
\end{eqnarray}

%and
%\begin{equation}\label{j-csillag-23}
%F^*(x,\alpha)=F(x,J^*(x,\alpha)). NEM KELL!!! \end{equation}

 Let $u:M\to
\mathbb R$ be a differentiable function in the distributional sense.
The {\it gradient} of $u$ is defined by
\begin{equation}\label{grad-deriv}
 \boldsymbol{\nabla}_F u(x)=J^*(x,Du(x)),
\end{equation}
where $J^*:T^*M\to TM$ is the {Legendre transform}
\begin{equation}\label{j-csillag}
    J^*(x,\alpha)=\sum_{i=1}^n\frac{\partial}{\partial
\alpha_i}\left(\frac{1}{2}F^{*2}(x,\alpha)\right)\frac{\partial}{\partial
x^i},
\end{equation} and $Du(x)\in T_x^*M$ denotes the (distributional) {\it
derivative} of $u$ at $x\in M.$ In local coordinates, one has
\begin{equation}\label{derivalt-local}
    Du(x)=\sum_{i=1}^n \frac{\partial u}{\partial x^i}(x){\rm d}x^i,
\end{equation}
$$\boldsymbol{\nabla}_F u(x)=\sum_{i,j=1}^n g_{ij}^*(x,Du(x))\frac{\partial u}{\partial x^i}(x)\frac{\partial}{\partial x^j}.$$
In general,
 $u\mapsto\boldsymbol{\nabla}_F u $ is not linear. The mean value
 theorem implies that
 \begin{equation}\label{novekedes}
    (Du(x)-Dv(x))(\boldsymbol{\nabla}_F u(x)-\boldsymbol{\nabla}_F
    v(x))\geq l_F F^{*2}(Du(x)-Dv(x)),\ \forall x\in M,
 \end{equation}
where
\begin{eqnarray}\label{lin-konst}
\nonumber l_F&\stackrel{\rm def.}{=}& \inf_{x\in M}\inf_{y,v,w\in T_xM\setminus \{0\}}\frac{g_{(x,v)}(y,y)}{g_{(x,w)}(y,y)} \\
 &=&\inf_{x\in
 M}\left(\frac{1-\|\beta\|_h(x)}{1+\|\beta\|_h(x)}\right)^2.
 \nonumber
\end{eqnarray}
 If $x_0\in M$ is
 fixed, then due to \cite{Ohta-Sturm}, one has
{ \begin{equation}\label{tavolsag-derivalt}
    F^*(x,D d_F(x_0,x))=F(x,\boldsymbol{\nabla}_F d_F(x_0,x))=D d_F(x_0,x)(\boldsymbol{\nabla}_F d_F(x_0,x))=1\ {\rm for\ a.e.}\ x\in
    M.
\end{equation}}
Let $X$ be a vector field on $M$. In a local coordinate system
$(x^i)$,  the {\it divergence} is defined by
div$_F(X)=\frac{1}{\sigma_F}\frac{\partial}{\partial x^i}(\sigma_F
X^i),$ where
 $\sigma_F(x)=\frac{\omega_n}{{\rm Vol}(B_x(1))}$, $\omega_n$ and Vol$(B_x(1))$ being
 the Euclidean volumes of the unit ball $B^n(1)$ and of the unit tangent ball $B_x(1)=\{y=(y^i):F(x,y^i \partial/\partial x^i)<
 1\}$, respectively.
The {\it Finsler-Laplace operator}
 $$\boldsymbol{\Delta}_F u={\rm div}_F(\boldsymbol{\nabla}_F u)$$ acts on $W^{1,2}_{\rm
 loc}(M)$ and for every $v\in C_0^\infty(M)$,
\begin{equation}\label{Green}
\int_M v\boldsymbol{\Delta}_F u {\text d}V_F(x)=-\int_M
Dv(\boldsymbol{\nabla}_F u){\text d}V_F(x),
\end{equation}
see
 \cite{Ohta-Sturm} and  \cite{Shen-monograph}. In particular,
 if $\beta=0$ (thus, the Randers space $(M,F)$ reduces to the Riemannian manifold $(M,h)$) the Finsler-Laplace operator
becomes the usual Laplace-Beltrami operator
 $\Delta_hu.$

 By definition, the {\it reversibility
constant} associated with $F$ is given by the formula
\begin{eqnarray}\label{reverzibilis}
 r_F&\stackrel{\rm def.}{=}& \sup_{x\in M}r_F(x),
\end{eqnarray}
where
\begin{eqnarray}\label{rev-unif-11}
\nonumber r_F(x)&\stackrel{\rm def.}{=}& \sup_{\substack{ y \in T_x
M\setminus \{0\}}}
\frac{F(x,y)}{F(x,-y)} \\
 &=&\frac{1+\|\beta\|_h(x)}{1-\|\beta\|_h(x)},
\end{eqnarray}
see Rademacher \cite{Rademacher}. It is clear that $r_{F}\geq 1$
(possibly, $r_{F}=+\infty$) and $r_{F}= 1$ if and only if $(M,F)$ is
reversible, i.e., $\beta=0$.

\subsection{The Funk-type metric $F_a$ on $B^n(1)$.} In this
subsection, we explicitly compute the objects introduced in the
previous subsection for the Funk-type metric $F_a$ on $B^n(1)$,
$a\in [0,1]$, given in (\ref{F-a-metrika}). Therefore, one has
$M=B^n(1)$ and the Randers metric $F_a$ is coming from the Klein
metric $h_K$,
$$(h_K)_x(y,y)=\frac{\sqrt{|y|^2-(|x|^2|y|^2-\langle
x,y\rangle^2)}}{1-|x|^2}$$ and by the 1-form
$$\beta_x=a\frac{ x}{1-|x|^2}.$$
It is clear that
$$(h_K)_{ij}=\frac{\delta_{ij}}{1-|x|^2}+\frac{x_ix_j}{(1-|x|^2)^2},\ i,j\in\{1,...,n\},$$
and according to Cheng and Shen \cite[Lemma 1.1.1]{Cheng-Shen}, we
have that $h_K^*=(h_K)^{-1}$ where the elements of the matrix are
given by
$$h_K^{ij}=(1-|x|^2)(\delta_{ij}-x_ix_j),\
i,j\in\{1,...,n\}.$$ Consequently,
\begin{equation}\label{beta-kiszamolva}
    \|\beta\|_{h_K}(x)=\sqrt{h_K^{ij}(\beta_x^i,\beta_x^j)}=a|x|.
\end{equation}
Therefore, the reversibility constant associated with $F_a$ on
$B^n(1)$ is given by
\[ r_{F_a}=  \left\{ \begin{array}{lll}
\frac{1+a}{1-a} &\mbox{if} &  a\in [0,1); \\
% u\geq 0 &\mbox{in} &   \Omega;\\
+\infty &\mbox{if} &  a= 1.
 \end{array}\right. \]
Due to (\ref{randers-volume}), we have
\begin{equation}\label{volume-form-hasonlitas}
    {\text d}V_{F_a}(x)=\left(1-a^2|x|^2\right)^\frac{n+1}{2}{\text
d}V_{h_K}(x),
\end{equation}
where the Klein volume form is given by $${\text
d}V_{h_K}(x)=\frac{1}{(1-|x|^2)^\frac{n+1}{2}}{\text d}x.$$ Finally,
the polar transform of $F_a$ is
\begin{equation}\label{F-a-polar}
    F_a^*(x,y)=\frac{\sqrt{(1-|x|^2)(1-a^2|x|^2)|y|^2-(1-a^2)(1-|x|^2)\langle x,y\rangle^2}-a(1-|x|^2)\langle
    x,y\rangle}{1-a^2|x|^2}.
\end{equation}
It is clear that $F_a^{**}=F_a$ and $r_{F_a^*}=r_{F_a}$.

\section{The Sobolev space on $(B^n(1),F_a)$: Proof of Theorem \ref{theorem-nem-vektor-ter}}

{\it Proof of Theorem \ref{theorem-nem-vektor-ter}.}
"(ii)$\Rightarrow$(i)." Let $a\in [0,1).$  Due to the convexity of
$F_a^{*2}$, if $u,v\in W_0^{1,2,a}(B^n(1))$ then $u+v\in
W_0^{1,2,a}(B^n(1)).$ Since $r_{F_a^*}=\frac{1+a}{1-a}$ is finite,
one also has that $c u\in W_0^{1,2,a}(B^n(1))$ for every $c\in
\mathbb R$ and $u\in W_0^{1,2,a}(B^n(1))$. The other properties of a
vector space are trivially verified.

"(i)$\Rightarrow$(ii)." The space $W_{0}^{1,2,a}(B^n(1))$ is assumed
to be a vector space over $\mathbb R$; by contradiction, we also
assume that one may have $a=1.$ In this case, $F_a$ is precisely the
Funk metric
$$F_1(x,y)=\frac{\sqrt{|y|^2-(|x|^2|y|^2-\langle
x,y\rangle^2)}}{1-|x|^2}+\frac{\langle x,y\rangle}{1-|x|^2},\ x\in
B^n(1),\ y\in \mathbb R^n.$$ Note that  the metric $F_1$ can be
obtained by
$$\left|x+\frac{y}{F_1(x,y)}\right|=1,$$
 while the distance function associated
to $F_1$ is given by
$$d_{F_1}(x_1,x_2)=\ln\frac{\sqrt{|x_1-x_2|^2-(|x_1|^2|x_2|^2-\langle x_1,x_2\rangle^2)}-\langle x_1,x_2-x_1\rangle}{\sqrt{|x_1-x_2|^2-(|x_1|^2|x_2|^2-\langle x_1,x_2\rangle^2)}-\langle x_2,x_2-x_1\rangle},\ x_1,x_2\in B^n(1),$$
see Shen \cite[p.141 and p.4]{Shen}. In particular,
$$d_{F_1}(0,x)=-\ln(1-|x|),\ x\in B^n(1).$$
First, by (\ref{tavolsag-derivalt}) or direct checking via
(\ref{F-a-polar}), we have that
\begin{equation}\label{alabb-kell---}
    F_1^*(x,Dd_{F_1}(0,x))=1.
\end{equation}
Second, a direct computation and (\ref{F-a-polar}) shows that
\begin{equation}\label{alabb-kell---2}
F_1^*(x,-Dd_{F_1}({0},x))=\frac{1+|x|}{1-|x|}.
\end{equation}

Let $u:B^n(1)\to \mathbb R$ be defined by
$u(x)=-\sqrt{1-|x|}=-e^{-\frac{d_{F_1}(0,x)}{2}}.$ It is clear that
$u\in W^{1,2}_{\rm loc}(B^n(1))$. First, since ${\rm
d}V_{F_1}(x)={\rm d}x$, we have
$$\int_{B^n(1)}u^2(x){\rm
d}V_{F_1}(x)=\frac{\omega_n}{n+1}.$$  On one hand, since
$Du(x)=\frac{1}{2}e^{-\frac{d_{F_1}({0},x)}{2}}D d_{F_1}({0},x),$ by
 (\ref{alabb-kell---}) it yields
\begin{eqnarray*}
  C_1 &:=&  \int_{B^n(1)}
F_1^{*2}(x,Du(x)){\rm d}V_{F_1}(x)=\frac{1}{4}\int_{B^n(1)}
{(1-|x|)}{\rm d}x\\
   &=&\frac{\omega_n}{4(n+1)}.
\end{eqnarray*}
Therefore, $\|u\|_{W^{1,2,1}}^2=\frac{5\omega_n}{4(n+1)},$ so $u\in
W_0^{1,2,1}(B^n(1))$.

On the other hand,  relation (\ref{alabb-kell---2})  implies that
\begin{eqnarray*}
  C_2 &:=&  \int_{B^n(1)}
F_1^{*2}(x,-Du(x)){\rm d}V_{F_1}(x)=\frac{1}{4}\int_{B^n(1)}
\frac{(1+|x|)^2}{1-|x|}{\rm d}x\\
   &=& +\infty,
\end{eqnarray*}
i.e., $-u\notin W_0^{1,2,1}(B^n(1))$, contradicting our initial
assumption. \hfill $\square$

\begin{remark}\rm \label{remark=1}
 Let $a\in [0,1).$ For every $x\in B^n(1)$, one has $0<1-a^2\leq
1-a^2|x|^2\leq 1$; thus, the volume forms ${\text d}V_{F_a}(x)$ and
${\text d}V_{h_K}(x)$ generate equivalent measures. Moreover, one
also has
\begin{equation}\label{hasonlitas}
    \frac{1}{(1+a)^2}h_K^*(y,y)\leq F_a^{*2}(x,y)\leq
\frac{1}{(1-a)^2}h_K^*(y,y),\ x\in B^n(1), \ y\in \mathbb R^n.
\end{equation}
Consequently,
$$\frac{(1-a^2)^\frac{n+1}{4}}{1+a}\|u\|_{H_1^2}\leq \|u\|_{W^{1,2,a}}\leq \frac{1}{1-a}\|u\|_{H_1^2},\ u\in C_0^\infty(B^n(1)).$$
In particular, the topologies generated by the objects
$(W_0^{1,2,a}(B^n(1)),\|\cdot\|_{W^{1,2,a}})$ and
$(H_1^2(B^n(1)),\|\cdot\|_{H_1^2})$ are equivalent whenever $a\in
[0,1).$ Moreover, a result of Federer and Fleming \cite{FF} for the
Klein ball model $(B^n(1),F_0)$ states that
\begin{equation}\label{Fed-fle-egy}
    \int_{B^n(1)}u^2(x){\rm d}V_{h_K}(x)\leq
\frac{4}{(n-1)^2}\int_{B^n(1)}h_K^*(Du(x),Du(x)){\rm d}V_{h_K}(x),\
\forall C_0^\infty(B^n(1)).
\end{equation}
 Therefore, the norm $\|\cdot \|_{H_1^2}$
and the 'gradient' norm over the Klein metric model given by
$$u\mapsto \|u\|_{K}=\left(\int_{B^n(1)}h_K^*(Du(x),Du(x)){\rm d}V_{h_K}(x)\right)^\frac{1}{2}$$
are also equivalent, i.e.,
\begin{equation}\label{ujabb-ekviv}
\|u\|_K\leq \|u\|_{H_1^2}\leq
\left(1+\frac{4}{(n-1)^2}\right)^\frac{1}{2}\|u\|_K.
\end{equation}
\end{remark}

\begin{remark}\label{remark-lp}\rm
The space $W_0^{1,2,1}(B^n(1))$ is a closed and convex cone in
$L_1^2(B^n(1)),$ where $L_a^p(B^n(1))$ denotes the usual class of
measurable functions $u:B^n(1)\to \mathbb R$ such that
$$\|u\|_{L_a^p}=\left(\int_{B^n(1)}|u(x)|^p{\text d}V_{F_a}(x)\right)^\frac{1}{p}<\infty$$
whenever $1\leq p<\infty.$ $L_0^p$ will be denoted in the usual way
by $L^p.$ $L^\infty(B^n(1))$ denotes the class of essentially
bounded functions on $B^n(1)$ with the usual sup-norm
$\|\cdot\|_{L^\infty}.$
\end{remark}

\section{Problem $({\mathcal P}_{\lambda})$: Proof of Theorem \ref{theorem-sublinear}}

Due to (g1), it follows that $g(0)=0$. Therefore,  one can extend
the function $g$ to
 $\mathbb R$ by $g(s)=0$ for $s\leq 0;$ this extension will be
 considered throughout of this section.

 An element $u\in W_0^{1,2,a}(B^n(1))$ is a weak
solution of problem $(\mathcal P_\lambda)$ if $u(x)\to 0$ as $|x|\to
1$ and
\begin{equation}\label{weak-solution-def}
    \int_{B^n(1)}Dv(\boldsymbol{\nabla}_{F_a} u){\text
    d}V_{F_a}(x)=\lambda\int_{B^n(1)}\kappa(x)g(u(x))v(x){\text
    d}V_{F_a}(x),\ \forall v\in C_0^\infty(B^n(1)).
\end{equation}

{\it Proof of Theorem \ref{theorem-sublinear}.} (i) Let $u\in
W_0^{1,2,a}(B^n(1))$ be a weak solution of $(\mathcal P_\lambda)$.
By density reasons, in (\ref{weak-solution-def}) we may use $v=u$ as
a test-function, obtaining by (\ref{volume-form-hasonlitas}),
(\ref{Fed-fle-egy}) and (\ref{hasonlitas}) that
\begin{eqnarray*}
  \int_{B^n(1)} F_a^{*2}(x,Du(x)){\rm
d}V_{F_a}(x) &=& \int_{B^n(1)}Du(\boldsymbol{\nabla}_{F_a} u){\text
    d}V_{F_a}(x) \\
   &=& \lambda\int_{B^n(1)}\kappa(x)g(u(x))u(x){\text
    d}V_{F_a}(x) \\
   &\leq& \lambda c_g\|\kappa\|_{L^\infty}\int_{B^n(1)}u^2(x){\text
    d}V_{h_K}(x) \\
   &\leq& \frac{4\lambda c_g\|\kappa\|_{L^\infty}}{(n-1)^2}\int_{B^n(1)}h_K^*(Du(x),Du(x)){\text
    d}V_{h_K}(x) \\
    &\leq& \frac{4\lambda c_g\|\kappa\|_{L^\infty}(1+a)^2}{(n-1)^2(1-a^2)^\frac{n+1}{2}}\int_{B^n(1)}F_a^{*2}(x,Du(x)){\text
    d}V_{F_a}(x). \\
\end{eqnarray*}
Consequently, if $0\leq \lambda
<c_g^{-1}\|\kappa\|_{L^\infty}^{-1}\frac{(n-1)^2(1-a^2)^\frac{n+1}{2}}{4(1+a)^2},$
$u$ is necessarily 0.

(ii) The proof is divided into
 several steps.

{\sc Step 1.} (Variational setting) Due to Remark \ref{remark=1}, we
may consider the energy functional $\mathcal J_{\lambda}:
H_1^2(B^n(1))\to \mathbb R$ associated with $(\mathcal
P_{\lambda})$, i.e.,
$$\mathcal J_{\lambda}(u)=\frac{1}{2}\mathcal E(u)-\lambda\mathcal G(u),$$
    where
    $$\mathcal E(u)=\int_{B^n(1)} F_a^{*2}(x,Du(x)){\rm d}V_{F_a}(x)\ {\rm and}\ \mathcal G(u)=\int_{B^n(1)}\kappa(x)G(u(x)){\text
    d}V_{F_a}(x).$$
On account of (g1), the functional $\mathcal J_{\lambda}$ is
well-defined and of class $C^1$; moreover, by (\ref{j-csillag}), we
have that
$$\mathcal J_{\lambda}'(u)(v)=\int_{B^n(1)} \left[Dv(\boldsymbol{\nabla}_{F_a}u)(x)-\lambda\kappa(x)g(u(x))v(x)\right]{\rm d}V_{F_a}(x).$$
In particular,  $\mathcal J_{\lambda}'(u)=0$ if and only if
(\ref{weak-solution-def}) holds.

{\sc Step 2.} (Symmetrization) Although $H_1^2(B^n(1))$ can be
embedded into the Lebesgue space $L^p(B^n(1))$, $p\in [2,2^*)$, see
Hebey \cite[Proposition 3.7]{Hebey}, this embedding is not compact.
Thus, we consider the space of radially symmetric functions in
$H_1^2(B^n(1))$, i.e.,
$$H_r(B^n(1))=\{u\in H_1^2(B^n(1)):u(x)=u(|x|)\}.$$
By using a Strauss-type inequality,  Bhakta and Sandeep
\cite{B-Sandeep} proved that the embedding
$H_r(B^n(1))\hookrightarrow L^p(B^n(1))$ is compact for every $p\in
(2,2^*)$. [Note that in \cite{B-Sandeep} the Poincar\'e ball model
is used which is conformally equivalent to the Klein ball model.]
Moreover, for every $u\in H_r(B^n(1))$, the Strauss-estimate shows
that $u(x)\to 0$ as $|x|\to 1.$

 If we introduce the action of the orthogonal group $O(n)$
on $H_1^2(B^n(1))$ in the usual manner, i.e., $$(\tau
u)(x)=u(\tau^{-1}x),\ u\in H_1^2(B^n(1)),\ \tau\in O(n),\ x\in
B^n(1),$$ then the fixed point set of $O(n)$ on $H_1^2(B^n(1))$ is
precisely the space $H_r(B^n(1))$. Moreover, by using
(\ref{F-a-polar}), let us observe that $F_a^*$ is $O(n)-$invariant,
i.e.,
\begin{equation}\label{f-a-invarians}
    F_a^*(\tau x,\tau y)=F_a^*(x,y),\ \forall \tau\in O(n),\ x\in B^n(1),\ y\in
\mathbb R^n.
\end{equation}
Therefore, since $D(\tau u)(x)=(\tau ^{-1})^tDu(\tau ^{-1}x)=\tau
Du(\tau ^{-1}x)$, where $\cdot^t$ denotes the transpose of a matrix,
we have for every $ \tau \in O(n)$ and $u\in H_1^2(B^n(1))$ that
\begin{eqnarray*}
  \mathcal E(\tau u) &=& \int_{B^n(1)} F_a^{*2}(x,D(\tau u)(x)){\rm d}V_{F_a}(x) \\
   &=& \int_{B^n(1)} F_a^{*2}(x,\tau Du(\tau ^{-1}x)){\rm d}V_{F_a}(x)\ \ \ \ \ ({\rm change\ of\ var.}\ \tau ^{-1}x=z) \\
   &=& \int_{B^n(1)} F_a^{*2}(\tau z,\tau Du(z)){\rm
   d}V_{F_a}(\tau z)\ \ \ \ \ ({\rm see\ (\ref{f-a-invarians})\ and}\ {\rm
   d}V_{F_a}(\tau z)={\rm
   d}V_{F_a}(z)) \\
   &=& \int_{B^n(1)} F_a^{*2}(z,Du(z)){\rm
   d}V_{F_a}(z)\\
   &=& \mathcal E(u),
\end{eqnarray*}
i.e., $\mathcal E$ is  $O(n)-$invariant. Similar reasoning as above
shows that $\mathcal G$ is also $O(n)-$invariant and  $O(n)$ act
isometrically on $H_1^2(B^n(1))$, i.e.,
$$\mathcal G(\tau u)=\mathcal G(gu)\ {\rm and}\ \|\tau u\|_{H_1^2}=\|u\|_{H_1^2},\ \forall \tau \in O(n),\ u\in H_1^2(B^n(1)).$$
By the above properties it follows that $\mathcal J_\lambda$ is
$O(n)-$invariant. Therefore, the principle of symmetric criticality
of Palais (see Krist\'aly, R\u adulescu and Varga \cite[Theorem
1.50]{KRV-book}) implies that the critical points of
$$\mathcal R_\lambda=\mathcal J_\lambda|_{H_r(B^n(1))}$$ are also
critical points for the original functional $\mathcal J_\lambda$. In
addition, since $u(x)\to 0$ as $|x|\to 1$ for every $u\in
H_r(B^n(1))$, we conclude that it is enough to guarantee critical
points for the functional $\mathcal R_\lambda$ in order to find
radially symmetric, weak solutions for problem $(\mathcal
P_\lambda)$.

For simplicity, let $\mathcal E_r$ and $\mathcal G_r$ be the
restrictions of $\mathcal E$ and $\mathcal G$ to $H_r(B^n(1))$,
respectively.  In the sequel, we shall show that there are at least
two critical points for $\mathcal R_\lambda$ whenever $\lambda$
belongs to a suitable interval.

{\sc Step 3.} (Subquadraticity of $\mathcal G_r$) We claim that
\begin{equation}\label{subliner-sobolev}
    \lim_{\substack{u\in H_r(B^n(1)) \\
\|u\|_{H_1^2}\to 0}}\frac{\mathcal G_r(u)}{\|u\|_{H_1^2}^2}=\lim_{\substack{u\in H_r(B^n(1)) \\
\|u\|_{H_1^2}\to \infty}}\frac{\mathcal G_r(u)}{\|u\|_{H_1^2}^2}=0.
\end{equation}
 Due to (g1), for every $\varepsilon>0$ there exists
$\delta_\varepsilon\in (0,1)$ such that
\begin{equation}\label{becsles-11}
    0\leq |g(s)|\leq \frac{\varepsilon}{\|\kappa\|_{L^\infty}}|s|\  {\rm
for\ all}\ |s|\leq \delta_\varepsilon\ {\rm and}\ |s|\geq
\delta_\varepsilon^{-1},
\end{equation}
%where $s_+=\max(s,0).$
Fix $p\in (2,2^*);$ clearly, the function $s\mapsto
\frac{g(s)}{s^{p-1}}$ is bounded on
$[\delta_\varepsilon,\delta_\varepsilon^{-1}]$. Therefore, for some
$m_\varepsilon>0,$ one has that
\begin{equation}\label{f-nov-becs-gyenge}
0\leq |g(s)|\leq \frac{\varepsilon}{\|\kappa\|_{L^\infty}}|s| +
m_\varepsilon |s|^{p-1}\ \ {\rm for\ all}\ s\in \mathbb R.
\end{equation}
 Thus, for
every $u\in H_r(B^n(1)),$ it yields that
\begin{eqnarray*}
 0\leq |\mathcal G_r(u)| &\leq& \int_{B^n(1)}\kappa(x)|G(u(x))|{\text
    d}V_{F_a}(x)  \\
   &\leq& \int_{B^n(1)}\kappa(x)\left[\frac{\varepsilon}{2\|\kappa\|_{L^\infty}}u(x)^2 + \frac{m_\varepsilon}{p} |u(x)|^{p}\right]{\text
    d}V_{h_K}(x)  \\
   &\leq& \int_{B^n(1)} \left[\frac{\varepsilon}{2}u(x)^2 + \frac{m_\varepsilon}{p} \kappa(x)|u(x)|^{p}\right]{\text
    d}V_{h_K}(x)  \\
   &\leq& \frac{\varepsilon}{2}\|u\|_{H_1^2}^2+ \frac{m_\varepsilon}{p}
   \|\kappa\|_{L^\infty}S_{p}^{p}\|u\|_{H_1^2}^{p},
\end{eqnarray*}
where $S_p>0$ is the best embedding constant in
$H_r(B^n(1))\hookrightarrow L^p(B^n(1)).$ Thus, for every $u\in
H_r(B^n(1))\setminus\{ 0\},$
\begin{eqnarray*}
  0\leq \frac{|\mathcal G_r(u)|}{\|u\|_{H_1^2}^2} &\leq& \frac{\varepsilon}{2}+ \frac{m_\varepsilon}{p}
  \|\kappa\|_{L^\infty}S_{p}^{p}%\left(1+\frac{4}{(n-1)^2}\right)^\frac{p}{2}
  \|u\|_{H_1^2}^{p-2}.
\end{eqnarray*}
Since $p>2$ and $\varepsilon>0$ is arbitrarily small,  the first
limit in (\ref{subliner-sobolev}) follows once $\|u\|_{H_1^2}\to 0$
in $H_r(B^n(1))$.

Let $q\in (1,2).$ Since $g\in C(\mathbb R,\mathbb R)$, there also
exists a number $M_\varepsilon>0$ such that
$$0\leq \frac{|g(s)|}{s^{q-1}}\leq M_\varepsilon\ {\rm for\ all}\ s\in
[\delta_\varepsilon,\delta_\varepsilon^{-1}],$$ where
$\delta_\varepsilon\in (0,1)$ is from (\ref{becsles-11}). The latter
relation together with (\ref{becsles-11}) give that
%\begin{equation}\label{f-nov-becs-kemeny}
$$0\leq |g(s)|\leq \frac{\varepsilon}{\|\kappa\|_{L^\infty}}|s| + M_\varepsilon |s|^{q-1}\ \ {\rm for\ all}\ s\in \mathbb R.$$
%\end{equation}
Similarly as above, it yields that
\begin{eqnarray}\label{eps-kkk}
  0\leq |\mathcal G_r(u)| &\leq &\frac{\varepsilon}{2}\|u\|_{H_1^2}^2+ \frac{M_\varepsilon}{q} \|\kappa\|_{L^\frac{2}{2-q}}\|u\|_{H_1^2}^{q}.
\end{eqnarray}
For every $u\in H_r(B^n(1))\setminus\{ 0\},$ we have that
\begin{eqnarray*}
  0\leq \frac{|\mathcal G_r(u)|}{\|u\|_{H_1^2}^2} &\leq& \frac{\varepsilon}{2}+ \frac{M_\varepsilon}{q}
  \|\kappa\|_{L^\frac{2}{2-q}}%\left(1+\frac{4}{(n-1)^2}\right)^\frac{r}{2}
  \|u\|_{H_1^2}^{q-2}.
\end{eqnarray*}
 Since $\varepsilon>0$ is arbitrary and $q\in (1,2)$, taking the limit $\|u\|_{H_1^2}\to
 \infty$ in $H_r(B^n(1))$,  we obtain the second relation in (\ref{subliner-sobolev}).

{\sc Step 4.} (Properties of $\mathcal R_\lambda$)  We are going to
prove that the functional $\mathcal R_\lambda$ is bounded from
below, coercive, and verifies the Palais-Smale condition on
$H_r(B^n(1))$ for every $\lambda \geq 0$. First, by (\ref{eps-kkk}),
it follows that
\begin{eqnarray*}
  \mathcal R_\lambda(u) &=& \frac{1}{2}\mathcal E_r(u)-\lambda\mathcal G_r(u) \\
   &\geq & \frac{(1-a^2)^\frac{n+1}{2}}{2(1+a)^2}\|u\|_K^2-\lambda\frac{\varepsilon}{2}\|u\|_{H_1^2}^2- \lambda\frac{M_\varepsilon}{r}
\|\kappa\|_{L^\frac{2}{2-r}}\|u\|_{H_1^2}^{r}.
\end{eqnarray*}
Choosing $\varepsilon>0$ sufficiently small, since
$\|\cdot\|_{H_1^2}$ and $\|\cdot\|_K$ are equivalent norms (see
(\ref{ujabb-ekviv})) and $r<2$, it follows that $\mathcal R_\lambda$
is bounded from below and coercive.

Now, let $\{u_k\}$ be a sequence in $H_r(B^n(1))$ such that
$\{\mathcal R_{\lambda}(u_k)\}$ is bounded and $\|\mathcal
R'_{\lambda}(u_k)\|_{*}\to 0.$ Since  $\mathcal R_{\lambda}$ is
coercive, the sequence $\{u_k\}$ is bounded in $H_r(B^n(1))$.
Therefore, up to a subsequence, we may suppose that $u_k\to u$
weakly in $H_r(B^n(1))$ and $u_k\to u$ strongly in $L^p(B^n(1))$ for
some $u\in H_r(B^n(1))$ and $p\in (2,2^*).$ In particular, we have
that
\begin{equation}\label{PS-1}
     \mathcal R_\lambda'(u)(u-u_k) \to 0\ {\rm and}\  R_\lambda'(u_k)(u-u_k) \to
     0\ {\rm as}\ k\to \infty.
\end{equation}
A direct computation gives that  $$
\int_{B^n(1)}(Du(x)-Du_k(x))(\boldsymbol{\nabla}_F
u(x)-\boldsymbol{\nabla}_F
    u_k(x)){\text
    d}V_{F_a}(x)= $$ $$=  \mathcal
R_\lambda'(u)(u-u_k)- R_\lambda'(u_k)(u-u_k)
+\lambda\int_{B^n(1)}\kappa(x)[g(u_k)-g(u)](u_k-u){\text
    d}V_{F_a}(x).
$$
By (\ref{PS-1}), the first two terms tend to zero. Moreover, due to
(\ref{f-nov-becs-gyenge}), it follows that
\begin{eqnarray*}
  T &:=& \int_{B^n(1)}\kappa(x)|g(u_k)-g(u)|\cdot|u_k-u|{\text
    d}V_{F_a}(x) \\
   &\leq & \int_{B^n(1)}\left(\varepsilon(|u_k|+|u|)+m_\varepsilon\|\kappa\|_{L^\infty}(|u_n|^{p-1}+|u|^{p-1})\right)|u_k-u|{\text
    d}V_{h_K}(x) \\
  &\leq &
  \varepsilon(\|u_k\|_{H_1^2}+\|u\|_{H_1^2})\|u_k-u\|_{H_1^2}\\&&+m_\varepsilon\|\kappa\|_{L^\infty} (\|u_k\|_{L^p}^{p-1}+\|u\|_{L^p}^{p-1})\|u_n-u\|_{L^p}.
\end{eqnarray*}
Since $\varepsilon>0$ is arbitrary small and $u_k\to u$ strongly in
$L^p(B^n(1))$, the last expression tends to zero.

 Moreover, relation
(\ref{novekedes}) implies that
\begin{eqnarray*}
  \mathcal E_r(u-u_k) &=& \int_{B^n(1)}F_a^{*2}(Du(x)-Du_k(x)){\text
    d}V_{F_a}(x) \\
   &\leq&\left(\frac{1+a}{1-a}\right)^2\int_{B^n(1)}(Du(x)-Du_k(x))(\boldsymbol{\nabla}_{F_a}
u(x)-\boldsymbol{\nabla}_{F_a}
    u_k(x)){\text
    d}V_{F_a}(x).
\end{eqnarray*}
Therefore, $\mathcal E_r(u-u_k)\to 0$ as $k\to \infty$, which means
in particular (see Remark \ref{remark=1}) that $\{u_k\}$ converges
strongly to $u$ in $H_r(B^n(1))$.

{\sc Step 5.} (First solution)
 On account of the assumption made on $\kappa$  and
(g2), there exists a truncation function $u_0\in
H_r(B^n(1))\setminus \{0\}$ such that $\mathcal G_r(u_0)>0.$ Thus,
we may define
$$\tilde \lambda=\inf_{ \substack{u\in H_r(B^n(1)) \\
\mathcal G_r(u)>0}}\frac{\mathcal E_r(u)}{2\mathcal G_r(u)}.$$ By
(\ref{subliner-sobolev}), we clearly have that $0<\tilde
\lambda<\infty.$ If we fix $\lambda>\tilde \lambda$, there exists
$\tilde u_\lambda\in H_r(B^n(1))$ with $\mathcal G_r(\tilde
u_\lambda)>0$ such that $\lambda>\frac{\mathcal E_r(\tilde
u_\lambda)}{2\mathcal G_r(\tilde u_\lambda)}\geq\tilde \lambda$.
Thus,
$$c_\lambda^1:=\inf_{H_r(B^n(1))}\mathcal R_\lambda\leq \mathcal R_\lambda(\tilde u_\lambda)=\frac{1}{2}\mathcal E_r(\tilde u_\lambda)-\lambda\mathcal G_r(\tilde u_\lambda)< 0.$$
Since $\mathcal R_\lambda$ is bounded from below and verifies the
Palais-Smale condition, the number $c_\lambda^1$ is a critical value
of $\mathcal R_\lambda$, i.e.,  there exists $u_\lambda^1\in
H_r(B^n(1))$ such that $\mathcal
R_\lambda(u_\lambda^1)=c_\lambda^1<0$ and $\mathcal
R_\lambda'(u_\lambda^1)=0.$ In particular, $u_\lambda^1\neq 0.$

{\sc Step 6.} (Second solution) Fix $\lambda>\tilde \lambda$.
Applying  (\ref{f-nov-becs-gyenge}) with the choice
$\varepsilon:=\frac{(1-a^2)^\frac{n+1}{2}}{2(\lambda+1)(1+a)^2}$, it
follows that
\begin{eqnarray*}
  \mathcal R_\lambda(u) &=& \frac{1}{2}\mathcal E_r(u)-\lambda\mathcal G_r(u) \\
   &\geq & \frac{(1-a^2)^\frac{n+1}{2}}{4(1+a)^2}\|u\|_K^2- \lambda\frac{m_\lambda}{p}
\|\kappa\|_{L^\infty}S_p^p\|u\|_{H_1^2}^{p},
\end{eqnarray*}
where $p\in (2,2^*)$ and $m_\lambda=m_\varepsilon>0$. Let
$$\rho_\lambda=\min\left\{\|\tilde u_\lambda\|_{H_1^2},\left(
\frac{(1-a^2)^\frac{n+1}{2}}{8\lambda\|\kappa\|_{L^\infty}S_p^p
m_\lambda(1+a)^2(1+\frac{4}{(n-1)^2})}\right)^\frac{1}{p-2}
\right\}.$$ The latter estimate and Step 5 shows that
$$\inf_{\|u\|_{H_1^2}=\rho_\lambda}\mathcal R_\lambda(u)=\eta_\lambda>0=\mathcal R_\lambda(0)> \mathcal R_\lambda(\tilde
u_\lambda),$$ i.e., the functional $\mathcal R_\lambda$ has the
standard mountain pass geometry. According to Step 4, one may apply
the mountain pass theorem, showing that there exists $u_\lambda^2\in
H_r(B^n(1))$ such that $\mathcal R_\lambda'(u_\lambda^2)=0$ and
$\mathcal R_\lambda(u_\lambda^2)=c_\lambda^2$, the number
$c_\lambda^2$ being characterized by
$$c_\lambda^2=\inf_{\gamma\in \Gamma}\max_{t\in [0,1]}\mathcal R_\lambda(\gamma(t)),$$
where $$\Gamma=\{\gamma\in C([0,1];H_r(B^n(1))):\gamma(0)=0,\
\gamma(1)=\tilde u_\lambda\}.$$ Since $c_\lambda^2\geq
\inf_{\|u\|_{H_1^2}=\rho_\lambda}\mathcal R_\lambda(u)>0$, it is
clear that $0\neq u_\lambda^2\neq u_\lambda^1.$ Since $g(s)=0$ for
every $s\leq 0,$ both solutions $u_\lambda^1$ and $u_\lambda^2$ are
non-negative, following from (\ref{weak-solution-def}).  This
concludes the proof. \hfill $\square$

\begin{remark}\rm (i) The case $a=1$ (Funk model) is not well
understood, since the set $W^{1,2,1}_0(B^n(1))$ is not a vector
space over $\mathbb R.$ However, we believe that variational
problems can also be treated within this context by using elements
from the theory of variational inequalities involving the indicator
function associated with the closed convex cone
$W^{1,2,1}_0(B^n(1))$ in $L^2_1(B^n(1))$, see \cite[Section
2]{KRV-book}.

(ii) For simplicity of the presentation, we considered elliptic
problems involving sublinear terms at infinity. The above
variational arguments seem to work also for elliptic problems
involving the Finsler-Laplace operator $\boldsymbol{\Delta}_{F_a}$,
$a\in [0,1)$, and superlinear or oscillatory nonlinear terms, see
e.g. Krist\'aly \cite{Kristaly-JDE}.
\end{remark}

\noindent {\bf Acknowledgment.} A. Krist\'aly is supported by a
grant of the Romanian National Authority for Scientific Research,
CNCS-UEFISCDI, "Symmetries in elliptic problems: Euclidean and
non-Euclidean techniques", project no. PN-II-ID-PCE-2011-3-0241.

%\section{Singular elliptic problems on hyperbolic spaces}

\end{document}